\documentclass[11pt]{article}
\usepackage{amssymb}
\input ss.macros
\newtheorem {thm}{Theorem}[section]
\newtheorem {lem}[thm]{Lemma}

\newtheorem {cor}[thm]{Corollary}

\def\Cox{\hfill \Box}

\def\ee{\epsilon}

\def\XX{{\bf X}}

\def\E{{\mathbb E}}
\def\P{{\mathbb P}}
\def\D{{\cal D}}
\def\diseq{\stackrel{\D}{=}}
\def\cE{{\cal E}}
\def\one{{\bf 1}}
\def\Var{{\rm Var}}

\def\br{\hfill \break}
\def\F{{\cal F}}
\def\|{\, | \, }

\begin{document}

\begin{titlepage}
\begin{center}
{\large \bf Cycles in random $k$-ary maps and the poor performance of 
random random number generation} 
\end{center}
\vspace{5ex}
\begin{flushright}
Robin Pemantle \footnote{Research supported in part by
National Science Foundation grant \# DMS 0103635}$^,$\footnote{The Ohio State
University, Department of Mathematics, 231 W. 18th Avenue, Columbus, OH
43210, pemantle@math.ohio-state.edu}
  ~\\
\end{flushright}

\vfill

{\bf ABSTRACT:} \br
Knuth~\cite{Kn} shows that iterations of a random function perform
poorly on average as a random number generator.  He proposes a
generalization in which the next value depends on two or more previous 
values.  This note demonstrates, via an analysis of the cycle length
of a random $k$-ary map, the equally poor performance
of a random instance in Knuth's more general model.  
\vfill



\end{titlepage}

\setcounter{equation}{0}
\section{Introduction} \label{ss:intro}

\subsection{Statement of problem}

In the introduction to his second volume, Knuth~\cite{Kn} discusses
the computer generation of pseudo-random numbers.  He gives several
cautionary tales about poor methods of generating these, including
a function whose description is so complicated that it mimics iterations
of a function chosen at random from all functions from 
$\{ 1 , \ldots , 10^{10} \}$ to itself.  The exercises
(see Exercises~11--15 on page 8 of~\cite{Kn}) then lead one through 
an analysis of a model where a function from $ [m] := \{ 1 , \ldots , m \}$ 
to itself is chosen uniformly at random.  The poor performance of this
pseudo-random number sequence is related to the cycle structure of
a random map and is well understood.  In particular, one may see
readily that the average length of the cycle of numbers produced from
a random seed is of order $\sqrt{m}$ and the cycle length from the
best seed is not much longer.  

Knuth then proposes the following generalization~\cite[Problem 19, 
page 9, labeled M48]{Kn}.  
A function $f$ is chosen uniformly from among the $m^{(m^k)}$ 
functions from $[m]^k$ to $[m]$.  Given an initial vector 
of values in $[m]^k$ for $(X_1 , \ldots , X_k)$, an infinite 
sequence of values is produced by the rule 
\begin{equation} \label{eq:recursion}
X_{n+k} = f(X_n , \ldots , X_{n+k-1}) \, .
\end{equation}
The problem is to determine the average length of the period of 
this eventually periodic sequence if the initial $k$ seeds are
chosen at random, and to answer as well some related questions: 
what is the chance that the eventual period has length 1, 
what is the average maximum cycle length over all seeds, what 
is the chance that there is no seed giving a cycle of length 1, 
and what is the average number of distinct eventual cycles as
the seed varies?  

\subsection{Heuristic}

A thumbnail computation shows that one might expect equally poor
performance from this multiple dependence model.  Let $W_n \in [m]^k$ 
denote $(X_n , \ldots , X_{n+k-1})$ and let $\mu < \tau$ be 
such that $W_{\tau - k + 1} = W_{\mu - k + 1}$ but the values of $W$ 
up to $W_{\tau - k}$ are distinct; thus the eventual period
is $\tau - \mu$ and the length of the sequence of $X$ values 
before repeating is $\tau - k$.  Although the values
$\{ W_n : n \geq 0 \}$ are no longer independent in the generalized model,
one may hope that they are nearly independent, so that the value
of the random quantity $\tau$ is well approximated
by the number $N$ of IID uniform draws from a population of size $M := m^k$
needed to obtain the first repeated value.  This is the classical
``birthday problem'' (see Example~(3d) on page~33 and the discussion
on page~49 of~\cite{Fe}).  It is known that the mean of $N$ in 
the birthday problem is of order $\sqrt{M}$ and more precisely
that $N^2 / 2M$ converges in distribution to a mean-one exponential as 
$M \to \infty$.  One would therefore expect (Theorem~\ref{th:2} below)
that $\tau^2 / 2 m^{k/2}$ converges to a mean-one exponential as well.

\subsection{Background}

The problem of random number generation is of coure fundamental
to the theory of computing.  There are many classes of problems,
including various factoring, counting and optimization problems,
for which randomized algorithms give solutions much faster (on
average) then any known non-randomized algorithm.  The use of
randomization in practice, is if anything, more widespread than
would be justified by theoretical results.  Monte Carlo methods
are ubiquitous in the areas of scientific computing, for example,
and banks use vast tables of pre-generated random numbers to price
derivative securities.

The need for sources of effectively random numbers has mushroomed 
with the explosion in computational ability.  Meanwhile, just as 
in the 1950's, the best sources are pseudo-random number generators, 
which are, for all their potential flaws, less prone to misbehave
than are physical sources of randomness~\cite{Ni}.

To construct objects with random or chaotic properties, it is often
easiest to pick them at random.  Expander graphs, for example, 
are easy to construct at random but difficult to construct
deterministically.  When it comes to pseudo-random number generation,
it is particularly appealing to use random or generic generators.
The problem of finding a class of generators, most of which 
generate good pseudo-random sequences, is therefore one of great
interest to theorists and practitioners alike.  

The random unary map is a natural starting point for exploration
of random pseudo-random sequences, but also underlies many other
phenomena in probability theory, ranging from random trees to
Brownian paths.  For this reason it has been studied in great
detail (see~\cite{Ko,AB,FO,AP1}) often with specific attention to
the time before repetition, as in~\cite{AP1}.  Its short cycle
time, hence unsuitability for random number generation, has
long been understood, as is evident from the discussion 
in all editions of~\cite{Kn}.  This has led to a huge industry
in random number generation (see, e.g.,~\cite{Ni}).  Classes 
of functions on $[m]$ or $[m]^k$ such as linear feedback registers 
or congruential sequences are sought whose periods are much 
longer than periods of functions chosen at random.  Interestingly,
the short cycle length can be a boon rather than always a drawback.
Pollard's rho-algorithm~\cite{Po} relies on the $\sqrt{m}$ cycle 
length to find any prime factor $p$ of $n$ in time $O(\sqrt{p})$;
Pollard's heuristic argument is in fact borne out empirically 
(see the discussion in~\cite[pages~466--471]{SF}).

Researchers studying random number generation via $k$-ary maps 
appear to have taken for granted that, as in the unary case, these
random maps cycle in a relatively short time.  Settling this,
however, was a problem stated already in the 1981 edition of~\cite{Kn},
given a rating of [M48], and left unsolved.  

\subsection{Results}

The purpose of this note is to show that the thumbnail computations
are correct.  All of the questions posed in Knuth may be correctly 
answered using the independence heuristic.  

In the case $k=2$ the following result holds.  Let $m$ be given
and let $\tau$, as defined above, be the least value for which
$W_{\tau - 1}$ repeats a previous value $W_j$ for $j < \tau - 1$.  

\begin{thm} \label{th:1}
As $m \to \infty$, the quantity $\frac{\tau^2}{2 m^2}$ converges
in distribution to an exponential of mean~1.  Furthermore, all
moments of $\frac{\tau^2}{2 m^2}$ converge to moments of the
exponential.  In particular,
\begin{eqnarray*}
\E  \; \tau & \sim & m \sqrt{\frac{\pi}{2}} \; ; \\
\Var (\tau) & \sim & \left ( 2 - \frac{\pi}{2} \right ) m^2 \; .
\end{eqnarray*}
\end{thm}

For general $k$ we will show:
\begin{thm} \label{th:2}
For any fixed $k$ and $x$, as $m \to \infty$, 
$$\P (\frac{\tau^2}{2 m^k} \geq x) \to \exp (-x) \, .$$
\end{thm}

Despite the fact that the arguments are straightforward, a careful
analysis may be justified for several reasons.  First, the question
has gone unanswered for sufficiently long (and not for lack of
interest) that the methods of analysis, though straightforward, must not 
be readily apparent.  It will therefore be useful to introduce to
the computer science community two techniques that are well known
to probabilists and theoretical statisticians.  In order to illustrate
the range of available techniques for this kind of analysis, two
different proofs will be presented.  

The first is a direct, combinatorial analysis and will be presented 
for the case $k=2$ (as is stated in Problem~16 to be the first 
interesting generalization), though it can easily be generalized 
to larger $k$.  It relies on the concept of {\em hazard rate}, 
well known in actuarial circles.  A brief introduction to this
concept is given in the next section.  The second analysis uses 
the Poisson approximation machinery of~\cite{AGG}, which relies 
on some technical lemmas of~\cite{BE} and concepts developed by 
Chen and Stein in the 1980's.  This method is discussed in
Section~3.  Although the Chen-Stein method is not elementary, 
the present application of this machinery is straightforward.  

A second reason for undertaking this analysis is re-inject 
questions about basic random models into the stream of scientific 
discussion. Often a result on a basic model 
will rekindle interest in simple variations and lead to a branch 
of research previously overlooked by a community that tends towards
depth-first research agendas.

Lastly, and perhaps most importantly, understanding
the behavior of iterations of a random $k$-ary map may prove
useful for other random models.  Just as short cycles of random 
unary maps have been used in Pollard's rho-algorithm and elsewhere,
it is not hard to imagine the random $k$-ary map underlying various 
other structures, algorithms and heuristics.

\setcounter{equation}{0}
\section{Time before repetition when $k=2$}

Given any sequence $\XX := \{ X_1 , X_2 , \ldots \}$ of values in $[m]$, 
we define the positive integer $\tau = \tau (\XX)$ as above to be 
minimal so that $W_{\tau - k + 1} = W_{\mu - k + 1}$ for some
$k \leq \mu < \tau$ (where $W_j$ are sub-words of length $k$ 
of the $\XX$ vector, as in the introduction).  When the $\XX$
vector is random, we let $\F_n$ denote the $\sigma$-field
$\sigma (X_1 , \ldots , X_n)$.  Compare the distributions of 
the vector $(\tau , X_1 , \ldots , X_\tau)$ under two different
measures for $\XX$: (a) when $\XX$ satisfies the 
recursion~(\ref{eq:recursion}) with $X_0 , \ldots , X_{k-1}$ 
IID uniform on $[m]$ and (b) when $\XX$ is an IID sequence 
of uniform draws from $[m]$.  Under both (a) and (b), the
conditional probability of $X_{n+1} = j$ given $\F_n$ is
$1/m$ as long as $\tau > n$.  The vector $(\tau , X_1 , \ldots , X_\tau)$ 
therefore has the same distribution under either law on $\XX$.
The main subject of our analysis it the distribution of $\tau$
and other quantities measurable with respect to $\F_\tau$.  We
will therefore assume throughout that $\XX$ is an infinite
IID uniform sequence.

Let $m$ be given.  In the remainder of this section, 
the dependence length, $k$, is fixed at two; thus $\tau$ 
is minimal so that $(X_{\tau-1} , X_\tau) = (X_{\mu-1} , X_\mu)$ 
for some $\mu < \tau$.  Arguments will be based on the hazard
rate principle, an informal statement of this is the following.

\begin{quote} {\small
The exponential lieftime of mean~1 is charactarized by the property 
that at any given time $t$, if you haven't died yet, then you have 
chance $dt$ of dying in the interval $(t , t + dt)$.  Now suppose
we observe an individual and determine at every time $t$ a chance
$f(t) \, dt$ for him to die in the next $dt$.  This will be by 
defintion a random variable measureable with respect to events up
to time $t$ (weight, smoking habits, and so on).  It is no longer
an easy matter to determine the unconditional lifetime distribution
of the population from typical information about the random function
$f$, but we can say somthing under a random time change.  If
the lifetime always has a conditional future density, then
the amount of hazard up to the time of death, defined by
$\int_0^T f(t) \, dt$ where $T$ is the death time, will always
be exponentially distributed with mean~1 (see~\cite[Prop.~3.28]{Je}).
This principle is often used for bounding the tails 
$\P (T > x)$, e.g., by $e^{-c} + \P (f(x) < c)$.
}\end{quote}

The proof of Theorem~\ref{th:1} is an elementary chain of asymptotic 
equalities.  Letting $\cE (1)$ denote an exponential of mean~1, 
we will show that
$$\cE (1) \diseq \mbox{ hazard } \approx \mbox{ linearized hazard }
   \approx \frac{\tau^2}{2 m^2} \, .$$
The first equality states that in the present discrete-time context, 
one may still find a time change under which the lifetime is an 
exponential.  The trick is to introduce auxilliary randomness to 
form the fractional part of the lifetime.  This is Lemma~\ref{lem:hazard} 
below, which is true for any stopping time.  The remaining asymptotic
inequalities then deal with approximations introduced by the discrete
time steps and by a small amount of unpredictability in the hazard rate
(cf.\ the methods of~\cite{Pe}).  
\begin{lem} \label{lem:hazard}
Let $\tau > 0$ be a stopping time on a probability space $(\Omega , \P)$ 
with respect to a filtration $\{ \F_n : n \geq 0 \}$ 
and let $h(k)$ be the random variable defined by 
$h(k) = - \log (1 - A_k)$ on the event that $\tau > k$ and 
arbitrarily otherwise, where
$$A_k := \P (\tau = k+1 \| \F_k) \, .$$
Let 
$$h = h^* + \sum_{j=0}^{\tau - 2} h(j)$$
where $h^*$ is a random variable, whose conditional distribution
given $h(\tau - 1) = x$ is a mean~1 exponential conditioned to be less
than $x$.  Suppose $\sum h(j) = \infty$ almost surely.
Then $h$ is distributed exactly as a mean~1 exponential.
\end{lem}

\noindent{\sc Proof:} Given $x > 0$, let $G_n$ be the event that
$\sum_{j=0}^{n-2} h(j) < x \leq \sum_{j=0}^{n-1} h(j)$.  Then
$G_n \in \F_{n-1}$ and
\begin{eqnarray*}
\P (h \geq x) & = & \sum_{n=1}^\infty \P ( h \geq x , G_n) \\
& = & \sum_{n=1}^\infty \E \P ( h \geq x , G_n \| \F_{n-1} ) \\
& = & \sum_{n=1}^\infty \E \one_{G_n} \P ( h \geq x \| \F_{n-1} ) \\
& = & \sum_{n=1}^\infty \E \one_{G_n} \left ( \prod_{j=0}^{n-2} e^{-h(j)}
   \right ) e^{-(x - \sum_{j=0}^{n-2} h(j))} \\
& = & \sum_{n=1}^\infty \E \one_{G_n} e^{-x} \\
& = & e^{-x} 
\end{eqnarray*}
since $\sum_j h(j) = \infty$ implies $\sum_n \one_{G_n} = 1$.   $\Cox$

The remainder of the proof of Theorem~\ref{th:1} involves combinatorial
specification of the hazard rate.  Apply the hazard rate lemma to the
quantity $\tau$ in the statement of Theorem~\ref{th:1}, resulting in 
quantities $h(k)$ and $h$ satisfying $h \diseq \cE (1)$.  Define $Y_k$
to be the number of $j < k$ for which $X_j = X_k$.  An easy lemma is:
\begin{lem} \label{lem:max}
As $m \to \infty$, $m^{-1/2} \max \{Y_k : k \leq \tau \} \to 0$ in 
probability. 
\end{lem}

\noindent{\sc Proof:} Keep a tally of how many times each value 
has been seen in the sequence $X_1 , X_2 , \ldots$.  Since these
are independent draws, it is evident that with probability 
$O(e^{-c_\ee \sqrt{m}})$ for some $c_\ee > 0$, no value 
is taken on $\ee \sqrt{m}$ times before every value is taken on 
$\ee \sqrt{m}/2$ times.  At a time $T_\ee$ when every value has been taken 
on $\ee \sqrt{m} / 2$ times, the hazard function $h$ is at least $c (\ee) m$, 
where $c (\ee)$ is a constant not depending on $m$.  It follows 
for fixed $\ee > 0$ that the probability of $\tau > T_\ee$ is 
exponentially small in $m$, and consequently that the probability of 
$\max \{ Y_k : k \leq \tau \}$ exceeding $\ee \sqrt{m}$ is at most
the sum of two probabilities that are exponentially small in 
$\sqrt{m}$, and hence that it tends to zero as $m \to \infty$.   $\Cox$

Recast the definition of $A_k$ in terms of $Y_k$,
$$A_k = \P (\tau = k + 1 \| \F_k ) = \frac{Y_k}{m} \, ,$$
to obtain the following immediate consequence.
\begin{cor} \label{cor:h}
For every $\ee > 0$ there is a $c_\ee > 0$ such that
$$\left | h - \sum_{j=0}^{\tau - 1} h(j) \right | \leq \ee m^{-1/2}$$
with probability at least $1 - c_\ee^{-1} (\exp c_\ee \sqrt{m})$.
\end{cor}

\noindent{\sc Proof:} By the definition of $h$,
\begin{equation} \label{eq:h minor}
0 \geq h - \sum_{j=0}^{\tau-1} h(j) \geq - h (\tau - 1) = \log 
   (1 - \frac{Y_{\tau - 1}}{m} ) \, .
\end{equation}
By Lemma~\ref{lem:max} this is at most $\ee m^{-1/2}$ except on a set
of measure tending to zero exponentially in $\sqrt{m}$ for each fixed
$\ee$.   $\Cox$

The cumulative linearized hazard rate
$$H(k) := \sum_{j=1}^k A_j$$
is close to $\sum_{j=0}^k h(j)$ but easier to work with.  We will see that
$$\sum_{j=0}^k h(j) \approx H_k \approx \frac{{k \choose 2}}{m^2} \, .$$
To quantify the last approximation, for $j \leq m$, let $T_k (j)$ be the number 
of $i \leq k$ for which $X_i = j$.  Then, counting pairs of occurrences 
of each value, an alternate definition of $H(k)$ is:
$$H(k) = \sum_{j=1}^m \frac{{T_k (j) \choose 2}}{m} \, .$$

\begin{lem} \label{lem:H}
If $k^2/m \to \infty$, then 
$$\frac{H_k}{{k \choose 2} / m^2} \to 1$$
in probability. 
\end{lem}

\noindent{\sc Proof:} Denote the first moment, second moment 
and variance of $H_k$ by $\mu_k , S_k$ and $V_k$ respectively.  
We may compute these as follows.
$\mu_k = E {T_k (1) \choose 2}$.  We compute $\mu_k$
as the expected number of pairs $(i,j)$ of indices 
at most $k$ for which $X_i = X_j = 1$.  Clearly then
$$\mu_k = \frac{{k \choose 2}}{m^2} \, .$$

Compute $m^2 S_k$ as $m E T_k (1)^2 + m(m-1) E T_k(1) T_k(2)$.
Counting ordered pairs of unordered pairs for which $X_u = X_v = 1$
and $X_w = X_x = 2$, we see that
$$E T_k(1) T_k(2) = \frac{{k \choose 2}{k -2 \choose 2}}{m^4} \, .$$
In a similar way, allowing for $(w,x)$ to have two, one or zero
elements in common with $(u,v)$, we get that
$$E T_k (1)^2 = \frac{{k \choose 2}}{m^2} + \frac{{k \choose 2} 2(k-2)}
   {m^3} + \frac{{k \choose 2}{k-2 \choose 2}}{m^4} \, .$$
Summing gives
$$S_k = {k \choose 2}{k-2 \choose 2} (\frac{m(m-1)}{m^6} + \frac{m}{m^6})
   + \frac{{k \choose 2}2(k-2)}{m^4} + \frac{{k \choose 2}}{m^3} .$$
Then
$$V_k = S_k - \mu_k^2 = \frac{{k \choose 2}}{m^3} - \frac{{k \choose 2}} 
   {m^4} $$
and 
$$\frac{V_k}{\mu_k^2} = \frac{m-1}{{k \choose 2}} \, .$$
The lemma now follows from Chebyshev's inequality.   $\Cox$

\noindent{\em Remark:} In order to prove convergence of all moments,
one must estimate $\E H(k)^p$ for integers $p > 2$.  There is an
expansion analogous to the equation $m^2 S_k = m \E T_k (1)^2
+ m(m-1) \E T_k (1) T_k (2)$.  Say that a descending vector 
of positive integers is a partition of $m$ if $\lambda = 
(\lambda_1 , \ldots , \lambda_{\# \lambda})$ and 
$\sum_{j=1}^{\# \lambda} \lambda_j = m$.  Let $T_k^\lambda$
denote the product $\prod_{j=1}^{\# \lambda} T_k (j)^{\lambda_j}$.  Then
$$m^p \E H(k)^p = \sum_{\lambda} m^{\# \lambda} (1 + O(m^{-1})) 
   \E T_k^\lambda$$
where the sum runs over partitions of $m$.  Here, the multiplier
$m^{\# \lambda} (1 + O(m^{-1}))$ gives the number of ways of choosing
distinct $j_1 , \ldots , j_{\# \lambda}$.  When $\lambda$ is the 
partition $(1 , \ldots , 1)$, the leading term of the sum is
$$m^p \frac{\prod_{j=0}^{p-1} {k - 2j \choose 2}}{m^{2p}}$$
leading to a contribution of $(1 + O(m^{-1})) (k^2 / (2m^2))^p$.
For any other $\lambda$, $\E T_k^\lambda$ is a sum of terms of the
form $O(k/m)^a$ with $1 \leq a \leq p$.  Each of these terms appears
in $\E H(k)^p$ with the multiplier $m^{\# \lambda - p}$, 
which is $O(m^{-1})$.  The total number of these terms is bounded,
so it follows that when $k/m > \ee$, 
\begin{equation} \label{eq:H}
\E H(k)^p = (1 + O(m^{-1})) \left ( \frac{k^2}{2 m^2} \right )^p \, .
\end{equation}
In other words, for $k/m \geq \ee$, $2 m^2 H_k / k^2$ converges to~1 
in each $L^p$ as $m \to \infty$, uniformly in $k$.  

\noindent{\sc Proof of Theorem}~\ref{th:1}: 
Convergence in distribution will follow from a comparison of 
$H(k)$ and $\sum_{j=0}^k h(j)$.  From the definitions,
\begin{eqnarray}
\sum_{j=0}^{k \wedge \tau - 1} h(j) & = & \sum_{j=1}^{k \wedge (\tau - 1)} 
   - \log (1 - A_k) \nonumber \\
& = & \sum_{j=1}^{k \wedge (\tau - 1)} A_k + O(A_k)^2 \nonumber \\
& = & H(k \wedge (\tau - 1)) \left ( 1 + O \left ( \max_{j \leq k \wedge 
   (\tau - 1)} A_j \right ) \right ) \nonumber \\
& = & H(k \wedge (\tau - 1)) \left ( 1 + O \left ( \max_{j \leq k \wedge 
   (\tau - 1)} \frac{Y_j}{m} \right ) \right ) \, . \label{eq:hH}
\end{eqnarray}
By Lemma~\ref{lem:max}, this shows that 
$\sum_{j=1}^{\tau - 1} h(j) / H(\tau - 1) \to 1$
in probability as $m \to \infty$.  Since $\tau^2 / m \to \infty$ in
probability as $m \to \infty$, Lemma~\ref{lem:H} may be applied to 
show that
\begin{equation} \label{eq:h}
\frac{2 n^2}{\tau^2} \sum_{j=0}^{\tau - 1} h(j) \to 1
\end{equation}
in probability as $m \to \infty$.  Together with Corollary~\ref{cor:h},
this implies that
$$\frac{2 m^2}{\tau^2} h \to 1$$
in probability as $m \to \infty$, and convergence in distribution 
of $\tau^2 / (2 m^2)$ to $\cE (1)$ then follows from the 
hazard rate lemma.

To extend this to convergence of higher integral moments, 
argue as follows.  We know that 
\begin{equation} \label{eq:higher1}
\cE (1) \diseq h \, .
\end{equation}
Let $|| \cdot ||_p$ denote the $L^p$ norm.  From Corollary~\ref{cor:h}
we see that 
\begin{equation} \label{eq:higher2}
|| \frac{h}{\sum_{j=0}^{\tau - 1} h(j)} ||_p \to 1
\end{equation}
as $m \to \infty$.  Let $G$ be the event that $\max \{ Y_k : k < \tau \}$
is at most $m^{-1/2}$.  It was shown in the proof of Lemma~\ref{lem:max}
that the probability of $G^c$ decays exponentially in $\sqrt{m}$.  
It was already shown in~(\ref{eq:hH}) that 
$$\left | \frac{\sum_{j=0}^{\tau - 1} h(j)}{H(\tau - 1)} \right |
   = 1 + O(m^{-1/2})$$
on $G$, which, together with the decay of $\P (G^c)$ faster than 
any polynomial, leads to
\begin{equation} \label{eq:higher3}
||\frac{\sum_{j=0}^{\tau - 1} h(j)}{H(\tau - 1)} ||_p \to 1
\end{equation}
as $m \to \infty$.  Finally, the estimate~(\ref{eq:H}) in the case
$k^2 / m \geq \ee$ together with convergence of $\tau^2 / m$ to $\infty$ 
in probability and monotonicity of $H(k)$ in $k$ imply that
\begin{equation} \label{eq:higher4}
||\frac{H(\tau - 1)}{\tau^2 / (2 m^2)} ||_p \to 1
\end{equation}
as $m \to \infty$.  The chain~(\ref{eq:higher1})--(\ref{eq:higher4})
of asymptotic equivalences in $L^p$ proves the last statement of the
theorem.   $\Cox$

\setcounter{equation}{0}
\section{Analysis of $X_{n+k} = f(X_n , \ldots , X_{n+k-1})$ for any
$k \geq 2$ via Poisson approximation}

In this section we will prove Theorem~\ref{th:2}.
Let $N := \lfloor \sqrt{2 m^k x} \rfloor$.  Then
$\tau^2 / (2 m^k) \geq x$ if and only if the values of $W_n$
for $0 \leq n \leq N$ are distinct.  Let $Z$ be the number
of pairs $(i,j)$ for which $0 \leq i < j \leq N$ and 
Theorem~\ref{th:2} is an immediate consequence of:
\begin{lem} \label{lem:pois}
The total variation distance between the law of $Z$ and
a Poisson of mean $x$ is $o(1)$ as $m \to \infty$.
\end{lem}

The proof of $Z$ is via the Chen-Stein Poisson Approximation method.
The classical Poisson approximation result~\cite[page~137]{Du}
says that if events are 
\begin{itemize}
\item independent,
\item each has probability at most $\ee$, which is going to zero,
\item and the sum of the probaiblities converges to $\lambda$,
\end{itemize}
then the number that occur converges in law to a Poisson of mean $\lambda$.
The Chen-Stein method is a means of weakening the independence 
assumption.  It suffices that most events be independent, and the rest
not too dependent.  Stein~\cite{St} developed an abstract framework for
quantifying such limit laws.  Later authors such as~\cite{AGG} provided 
useful hypotheses for getting numerical bounds on the 
distance, in total variation, to a Poisson distribution. 

\noindent{\sc Proof of Lemma}~\ref{lem:pois}:  
For the duration of this proof, $\alpha$
and $\alpha'$ will be shorthand for $(i,j)$ and $(i' , j')$
respectively.  Let $S$ denote the set of $\alpha$ for which
$0 \leq i < j \leq N$.  Let $G_\alpha$ denote the event that
$W_i = W_j$.  Define $p_\alpha = \P (G_\alpha)$ and
$p_{\alpha \alpha'} = \P (G_\alpha \cap G_{\alpha'})$.  

Let $B(\alpha)$ be the set of $\alpha'$ for which $|y-y'| < k$
for some $y \in \{ i , j \}$ and $y' \in \{ i' , j' \}$.
Note that for $\alpha' \notin B(\alpha)$, the event $G_{\alpha'}$
that $W_{i'} = W_{j'}$ is measurable with respect to 
$\{ X_s : |s - i| , |s - j| \geq k \}$. Therefore, 
\begin{equation} \label{eq:0}
G_\alpha \mbox{ is independent of } \sigma(G_{\alpha'} : 
   \alpha' \notin B(\alpha)) \, .
\end{equation}
Define
\begin{eqnarray}
b_1 & := & \sum_\alpha \sum_{\alpha' \in B(\alpha)} p_\alpha p_{\alpha'} 
   \label{eq:1} \, ; \\
b_2 & := & \sum_\alpha \sum_{\alpha \neq \alpha' \in B(\alpha)} p_{\alpha 
   \alpha'} \label{eq:2} \, .
\end{eqnarray}
The quantities $b_1$ and $b_2$ are quantities appearing under the
same name in~\cite[Theorem~1]{AGG}; their quantity $b_3$
is zero due to the independence relation~(\ref{eq:0}).  The conclusion
of~\cite[Theorem~1]{AGG} is that $|\P (Z=0) - \exp (\E Z)| < b_1 + b_2$.
It remains to identify $\E Z$ and to bound $b_1$ and $b_2$ from above.

Observe first the claim that for any $\alpha$, $p_\alpha = m^{-k}$.  
This is obvious for $|i - j| \geq k$.  But in fact for any
$i$ and $j$, $G_\alpha$ occurs if and only if $X_{i+r} = X_{j+r}$
for $0 \leq r < k$.  For $j = i+s$ with $0 < s < k$, the values
of $X_i , \ldots , X_{i+s-1}$ may be chosen arbitrarily, and there
will be precisely one set of values of $X_{i+s} , \ldots , X_{i+s+k-1}$
for which $G_\alpha$ occurs, proving the claim.  It follows that
\begin{equation} \label{eq:pois1}
\lambda := \E Z = \sum_\alpha p_\alpha = m^{-k} {N \choose 2} = 
   (1 + o(1)) x \, .
\end{equation}

Observe next that the cardinality of $B(\alpha)$ is at most 
$8 k N$, since the number of pairs $(i' , j')$ with $i'$ within
$k$ of $i$ is at most $2 k N$, and similarly for the other three
possibilities.  It follows immediately that 
\begin{equation} \label{eq:pois2}
b_1 \leq (\sum_\alpha p_\alpha) (8 k N) m^{-k}
   \leq (8 \sqrt{2} + o(1)) k x^{3/2} m^{-k/2} \, .
\end{equation}

Finally, we bound $b_2$ from above.  Let $B_0 (\alpha)$
denote the set of $\alpha'$ for which both of $i'$ and $j'$
are within $k$ of either $i$ or $j$.  Then $|B_0 (\alpha)|
< (4k)^2$.  

\noindent{\bf Claim:} for $\alpha \neq \alpha' \in B_0 (\alpha)$, 
$$p_{\alpha , \alpha'} \leq m^{-k-1} \, .$$
Assume without loss of generality that $j < j'$, since
the other cases, $j > j'$, $i < i'$ and $i > i'$ are similar.  Then 
$$p_{\alpha , \alpha'} \leq p_\alpha \P (G_{\alpha'} \|
   X_n : n < j') \leq m^{-k} m^{-1} \, ,$$
proving the claim.  

For $\alpha' \notin B_0 (\alpha)$, one
conditions on $\{ X_s : |s-i| < k \mbox{ or } |s-j| < k \}$ 
to see that $p_{\alpha \alpha'} = m^{-2k}$.  One then has
\begin{eqnarray} 
b_2 & = & \sum_\alpha \left [ \sum_{\alpha' \in B_0 (\alpha)} 
   p_{\alpha \alpha'}
   + \sum_\alpha \sum_{\alpha' \in B_0 (\alpha)} p_{\alpha \alpha'} \right ]
   \nonumber \\
& \leq & \sum_\alpha \left [ 16 k^2 m^{-k-1} + (8 k N) m^{-2k} \right ] 
   \nonumber \\
& \leq & 8 k^2 N^2 m^{-k-1} + 4 k N^{3/2} m^{-2k} \nonumber \\
& = & (1+o(1)) (16 k^2 \lambda m^{-1} + (4 2^{3/4} \lambda^{3/2} k m^{-k})
   \label{eq:pois3}
\end{eqnarray}
as $m \to \infty$.  Combining~(\ref{eq:pois1})~-~(\ref{eq:pois3})
establishes that $b_1 + b_2 = o(1)$ and $\E Z - x = o(1)$, which
completes the proof of Theorem~\ref{th:2}.   $\Cox$

\section{Further discussion}

Let $U$ and $\cE (1)$ be independent with $U$ uniform on $[0,1]$
and $\cE (1)$ exponential of mean 1.  The following extension 
of the distributional convergence results may be proved.
Recall that $\mu$ is the index for which $X_\mu , \ldots , X_{\tau - 1}$
is the first full period of the eventually periodic sequence of
pseudo-random numbers.  
\begin{thm} \label{th:4}
As $m \to \infty$, the pair $(\mu , \tau)$ converges in distribution
to $(U \cE (1) , \cE (1))$.
\end{thm}

Complete proof of the extensions in this section will not be given,
but the argument, along the lines of the first analysis, is as follows.
Fix an integer $r$ and break the hazard rate for the occurrence of $\tau$
into $r$ components.  The $j^{th}$ component at time $n$ is the 
hazard rate for the occurrence of $\tau = n+1$ and $(j-1)/r \leq 
\mu < j/r$.  A lemma analogous to Lemma~\ref{lem:H} shows that 
the $r$ hazards accumulate at asymptotically equal rates, and a lemma
analogous to the hazard rate lemma then shows the asymptotic uniform
distribution of $\mu / \tau$ over the $r$ bins given 
$\lfloor r \, \tau \rfloor$.  Sending $r$ to infinity completes
the argument.  

An analysis of the probability of landing in a cycle of length~1
is easiest along the lines of the Poisson approximation.  Indeed,
the number of occurrences of $W_n$ of the form $(j , \ldots , j)$ for
some $j \leq m$ by time $k$ is well approximated by a Poisson of mean
$k m^{1-k}$; the number of these followed by one more $j$ is then
nearly a Poisson of mean $k m^{-k}$.  Since $\tau$ is of order
$m^{k/2}$, one sees that the mean number of these occurrences by
time $\tau$ is $\Theta (m^{-k/2})$, so this gives the order of
magnitude of the chance of being caught in a cycle of length~1.
On the other hand, the probability that some seed results in a 
cycle of length~1 is the chance that one of the $m$ words
$(j , \ldots , j)$ maps to itself, which rapidly approaches
$1 - e^{-1}$ as $m \to \infty$.  

An upper bound on the maximum value of $\tau$ over all seeds is obtained
as follows.  In the spirit of Theorems~\ref{th:1} and~\ref{th:2},
the probability that $\tau^2 > 2 (1 + \ee) m^k (k \log m)$
can be shown to be close to $\exp (- (1 + \ee) k \log m)$.
Indeed, while Theorems~\ref{th:1} and~\ref{th:2}, as written,
compute $\P (\tau^2 > (2 m^k) x)$ only when $x$ is fixed, the
arguments are sufficient to handle poly-logarithmic growth of $x$,
that is $x \leq (\log m)^p$.  Specifically, the four chains
in the asymptotic equalities when $x$ grows at this rate are:
the exact equality $h \diseq \cE (1)$ as before; the difference 
between $h$ and $\sum_{j-0}^{\tau - 1} h(j)$ small in every $L^p$;
the linearization error in replacing $h$ by $H$ changes the
likelihood of exceeding a hazard of $x$ from $e^{-x}$ to
$e^{-x + o(x)}$, and the ratio between $H(k)$ and its deterministic
counterpart $k^2 / (2 m^2)$ is small as long as $x$ is not too small
(as before).  One may then extend the estimate to slowly growing $x$:
$$\P (\tau^2 > 2 (1 + \ee) k m^k \log m ) \sim m^{- (1+\ee) k} \, .$$
Since there are $m^k$ seeds, this gives
\begin{equation} \label{eq:log}
\P \left [ \tau^* > \sqrt{b \, k \, m^k \, \log m} \right ] \to 0
\end{equation}
for any $b > 2$, where $\tau^*$ is the supremum over seeds of the 
value of $\tau$ for a fixed random $f$.  

A interesting theoretical problem, perhaps do-able with only 
a modest amount of effort, would be to make this more
precise and give a sharper estimate.  For iterations of a 
unary map, the distribution of $\tau^*$ is known.  Its mean
was shown in~\cite[Theorem~7]{FO} to be asymptotic to an
explicit contant multiple of $\sqrt{m}$.  Later the scaling
limit of $m^{-1/2} \tau^*$ was shown to exist~\cite{AP1} 
and an explicit formula given~\cite[Theorem~1]{AP2}.

\noindent{\bf Problem:} {\em Let $\tau^* (m,k)$ be the supremum over all 
seeds of the value of $\tau$ for iterations of one random function
$f : [m]^k \to [m]$.  Show that $m^{-k/2} \tau^*$ converges
weakly to the same limit as described in~\cite{AP1,AP2}.  }

If this proves difficult, perhaps it could at least be shown
that $\tau^* = O(m^{k/2})$ in probability, that is, that the
extra factor of log under the radical in~(\ref{eq:log}) is
superfluous.

\noindent{\bf Acknowledgement:} The author would like to
thank Philippe Flajolet for bringing this problem to light at
the $8^{th}$ Analysis of Algorithms meeting in Strobl.

\end{document}